\documentclass[11pt]{amsart}
\usepackage[a4paper,margin=26mm]{geometry}
\usepackage{amsmath,amssymb,amsthm,mathtools}
\usepackage[T1]{fontenc}
\usepackage{lmodern,microtype}
\usepackage{color}
\usepackage[colorlinks=true,linkcolor=blue,citecolor=blue,urlcolor=blue]{hyperref}
\usepackage{enumitem}
\setlist{itemsep=3pt,topsep=5pt}
\allowdisplaybreaks[2]
\numberwithin{equation}{section}
\newtheorem{theorem}{Theorem}[section]
\newtheorem{lemma}[theorem]{Lemma}
\newtheorem{proposition}[theorem]{Proposition}
\theoremstyle{remark}\newtheorem{remark}[theorem]{Remark}
\newcommand{\R}{\mathbb R}

\newcommand{\dd}{\,\mathrm d}

\hypersetup{pdftitle={A direct level-set proof of the fractional heat Liouville theorem},pdfauthor={Yahong Guo, Yugao Ouyang, and Jiongduo Xie}}

\begin{document}
\title[Liouville theorems for the fractional heat equation]
{An Optimal Liouville Theorem 
for the Fully Fractional Heat Equation}

\author[Guo]{Yahong Guo}
\address{School of Mathematical Sciences\\
Shanghai Jiao Tong University\\
Shanghai, 200240, China} \email{yhguo@sjtu.edu.cn}
\author[Ouyang]{Yugao Ouyang}
\address{School of Mathematical Sciences\\
Shanghai Jiao Tong University\\
Shanghai, 200240, China}  \email{ouyang1929@sjtu.edu.cn}
\author[Xie]{Jiongduo Xie}
\address{School of Mathematical Sciences\\
Shanghai Jiao Tong University\\
Shanghai, 200240, China}
\email{jiongduoxie@outlook.com
}

\begin{abstract}
We prove that every nonnegative entire classical solution of the fully fractional heat equation $(\partial_t-\Delta)^s u=u^p$ in $\mathbb{R}^n\times\mathbb{R}$ is identically zero for $n\geq2$, $0<s<1$, and $1<p<(n+2s)/(n-2s)$. No global boundedness, decay, or symmetry assumptions are required for the solutions. The exponent range is sharp, since the critical equation admits positive stationary solutions. Our result provides a fractional counterpart of Quittner's optimal Liouville theorem for the classical semilinear heat equation, while our proof follows a new and simpler approach. For bounded solutions, we combine elliptic nonexistence and compactness to control superlevel set measures by the local energy dissipation of the weighted extension. A local energy inequality then gives nonexistence by iteration. Doubling, local regularity, and a nonlocal convergence theorem yield a universal pointwise bound in terms of the parabolic distance to the boundary, which implies our Liouville theorem for general solutions. 
\end{abstract}
\maketitle

\section{Introduction and main results}
We study nonnegative entire solutions of the Lane--Emden type master
equation
\begin{equation}\label{master-lane-emden}
 H^s u:=(\partial_t-\Delta)^s u=u^p
 \quad\text{in }\R^n\times\R,
 \qquad n\geq2,\quad 0<s<1,\quad p>1,
\end{equation}
where $H=\partial_t-\Delta$ is the classical heat operator. The fractional heat operator is nonlocal in both space and time: its value at $(x,t)$ depends on the solution at every spatial point and
at all earlier times. More precisely, it  is defined by
\begin{equation}\label{nonlocaloper}
 H^s u(x,t)
 =C_{n,s}\int_{-\infty}^{t}\int_{\R^n}
 \frac{u(x,t)-u(y,\tau)}{(t-\tau)^{n/2+1+s}}
 e^{-\frac{|x-y|^2}{4(t-\tau)}}\dd y\dd\tau,
 \qquad
 C_{n,s}=\frac{1}{(4\pi)^{n/2}|\Gamma(-s)|}.
\end{equation}
Here $\Gamma$ is the Gamma function, 
the integral in $y$ is taken in the Cauchy principal value sense. Equivalently, setting
\[
 S_r u(x,t)=(4\pi r)^{-n/2}
 \int_{\R^n}e^{-|z|^2/(4r)}u(x-z,t-r)\dd z,
 \qquad r>0,
\]
the changes of variables $r=t-\tau$ and $z=x-y$ give
\begin{equation}\label{eq:operator}
 H^s u(x,t)=\frac{1}{|\Gamma(-s)|}
 \int_0^\infty\bigl[u(x,t)-S_r u(x,t)\bigr]
 \frac{\dd r}{r^{1+s}}.
\end{equation}
When acting on functions independent of time, $H^s$ reduces to $(-\Delta)^s$, whereas on functions independent of space it becomes the left-sided Marchaud derivative of order $s$. Moreover, as $s\uparrow1$, the fractional heat operator $H^s$ converges
to the classical heat operator $H=\partial_t-\Delta$, in the sense that
$H^s\varphi\to H\varphi$ uniformly for each
$\varphi\in C_c^\infty(\mathbb R^n\times\mathbb R)$. (cf. \cite{FNW2023}).

We now specify the solution class. Let
\begin{equation}\label{eq:tail-space}
\begin{split}
 \mathcal L_s(\R^n\times\R)
 =\Biggl\{u\in L^1_{\mathrm{loc}}(\R^n\times\R):\
 &\int_{-\infty}^{t}\int_{\R^n}
 \frac{|u(y,\tau)|e^{-|x-y|^2/[4(t-\tau)]}}
 {1+(t-\tau)^{n/2+1+s}}\dd y\dd\tau<\infty\\[-2pt]
 &\hspace{15mm}\text{for every }(x,t)\in\R^n\times\R
 \Biggr\}.
\end{split}
\end{equation}
For an open set $\Omega\subset\R^n\times\R$, a function $u$, defined
in $\R^n\times\R$, is a classical solution of $H^s u=u^p$ in $\Omega$
if $u\in\mathcal L_s(\R^n\times\R)\cap C^{2s+\varepsilon,\,s+\varepsilon/2}_{x,t,\mathrm{loc}}(\Omega)$ for some $\varepsilon>0,$
and the equation holds pointwise in $\Omega$.
When $\Omega=\R^n\times\R$, we call $u$ an entire solution.
Here and below, $C^{\alpha,\alpha/2}_{x,t}$ denotes the usual parabolic
H\"older space.

The fractional heat operator $(\partial_t-\Delta)^s$ arises in the study
of coupled continuous-time random walks, which have applications in
anomalous diffusion~\cite{MBSBK,FNW} and finance~\cite{RSM}.
It also appears in
obstacle problems for semipermeable membranes~\cite{ACM} and, for $s=1/2$,
in a boundary formulation of biological invasion along roads~\cite{BRR,ST2017}.
In this case, \eqref{master-lane-emden} describes half-space heat flow
with boundary flux $u^p$, a power law for heat release from surface
oxidation in radiative ignition models~\cite{BGKKS}. Moreover, both qualitative and quantitative aspects of the fractional heat equation have been extensively investigated in recent years, as reflected in the growing body of literature \cite{CS2014,CG1,CG2,CGL2026,CGLO2026,CGM2025,FdP2024,ST2017}.


Our goal is to establish an optimal Liouville theorem in the Sobolev-subcritical range for
global nonnegative solutions of \eqref{master-lane-emden}. The classical elliptic counterpart
goes back to Gidas and Spruck~\cite{GS1981}, who used integral identities
to prove that, for $n\geq 3$ and $1<q<(n+2)/(n-2)$, the equation
$-\Delta u=u^q$ admits no positive solutions in $\mathbb R^n$.
Chen and Li~\cite{CL1991} later gave a different proof using the method
of moving planes. For the master equation, the corresponding Sobolev
exponent arises naturally from its stationary reduction:
if $u$ is independent of time, then $H^s u=u^p$ reduces to
$(-\Delta)^s u=u^p$ in $\mathbb R^n$. Let
\begin{equation}\label{eq:sobolev}
 p(n,s)=\frac{n+2s}{n-2s}.
\end{equation} For $1<p<p(n,s)$, this equation admits no positive solutions, 
see Chen, Li, and Li~\cite{CLL2017} .
For the classical heat equation $u_t-\Delta u=u^p$,
Quittner~\cite{Quittner2021} proved that there are no positive entire classical
solutions throughout the Sobolev-subcritical range $p>1$, $p<(n+2)/(n-2)$,
without any boundedness, symmetry, or decay assumption. 
Using the Giga--Kohn monotonicity formula in similarity variables,
Quittner iteratively refined the energy estimates. After rescaling,
these estimates yield a positive stationary limit, contradicting the
elliptic Liouville theorem of Gidas and Spruck~\cite{GS1981}. Quittner~\cite[Theorem~1]{Q} adapted this approach to the heat equation
in a half-space with the nonlinear boundary condition $-\partial_y U=U^p$.
Through the heat extension~\cite{NS2016,ST2017}, his result gives the
Liouville theorem for bounded positive solutions of
\eqref{master-lane-emden} when $s=\tfrac12$ and $1<p<(n+1)/(n-1)$.



Recently, several  works related to Liouville-type theorems have focused on the fully fractional heat equation \eqref{master-lane-emden} considered here. 
 Ferreira and de Pablo~\cite{FdP2024} studied \eqref{master-lane-emden}  in $\mathbb{R}^n\times (0,T)$ with memory and proved finite-time
blow-up for nontrivial nonnegative solutions when
$1<p\le p_F(n,s)$, where
\begin{equation}\label{eq:fujita}
 p_F(n,s)=1+\frac{2s}{n+2-2s}<p(n,s)
\end{equation}
is the Fujita exponent. Their comparison argument also rules out  nontrivial nonnegative bounded entire solutions in $\mathbb{R}^n\times\mathbb{R}$ in this range.
The present work will extend this nonexistence result to the full
Sobolev-subcritical range $1<p<p(n,s)$, without assuming global boundedness.
Sire, Wei, Wu, and Ye~\cite{SWWY} considered ancient weak solutions of $H^s u=|u|^{p-1}u$ in $\R^n\times(-\infty,0)$ satisfying the global Type-I bound
\begin{equation}\label{eq:typeI}
 \sup_{x\in\R^n,\ t<0}(-t)^\beta|u(x,t)|<\infty,
 \qquad \beta=\frac{s}{p-1}.
\end{equation}
For $1<p\leq p(n,s)$, they classify self-similar solutions as $u\equiv0$ or $u(x,t)=\pm\kappa(-t)^{-\beta}$, where $\kappa=\left(\Gamma(\beta+s)/\Gamma(\beta)\right)^{1/(p-1)}.$ Besides, they can
remove self-similarity under the additional condition
\[
 \limsup_{t\uparrow0}(-t)^\beta|u(0,t)|>0.
\]

The present paper concerns entire solutions without imposing global boundedness, decay, symmetry, self-similarity, or a Type-I bound. We prove that every
nonnegative classical entire solution vanishes when $1<p<p(n,s)$.
The exponent range is sharp since, at $p=p(n,s)$, the positive bubbles of the fractional
Lane--Emden equation~\cite{CLZ2017} give stationary solutions of
\eqref{master-lane-emden}.

Our first main result is the following Liouville theorem.
\begin{theorem}\label{main-thm-1}
Let $n\geq2$, $0<s<1$, and $1<p<p(n,s)$.
If $u$ is a nonnegative classical entire solution of
\eqref{master-lane-emden}, then $u\equiv0$.
\end{theorem}

The above theorem is a direct consequence of the following local a priori estimate for solutions in any open set. To state the universal estimate, we use the parabolic distance
\[
 d_{\mathrm p}\bigl((x,t),(y,\tau)\bigr)
 =\max\{|x-y|,|t-\tau|^{1/2}\},
 \qquad
 \operatorname{dist}_{\mathrm p}(P,E)
 =\inf_{Q\in E}d_{\mathrm p}(P,Q).
\]
We also write
\[
 Q_R(x_0,t_0)=B_R(x_0)\times(t_0-R^2,t_0+R^2),
 \qquad B_R=B_R(0),\qquad Q_R=Q_R(0,0),
\]
where $B_R(x_0)$ is the Euclidean ball in $\R^n$.

\begin{theorem}\label{main-thm-2}
Let $n\geq2$, $0<s<1$, and $1<p<p(n,s)$.
There exists $C=C(n,s,p)>0$ such that, for every nonempty proper
open set $\Omega\subsetneq\R^n\times\R$ and every classical solution
of $H^s u=u^p$ in $\Omega$ satisfying
$u\geq0$ in $\R^n\times\R$, one has
\begin{equation}\label{estimate-universe}
 u(P)\leq
 C\,\operatorname{dist}_{\mathrm p}(P,\Omega^c)^{-\frac{2s}{p-1}}
 \qquad\text{for all }P\in\Omega.
\end{equation}
\end{theorem}

The power in \eqref{estimate-universe} agrees with the scaling
\[
 u_\lambda(x,t)=\lambda^{\frac{2s}{p-1}}
 u(\lambda x,\lambda^2t),
 \qquad \lambda>0,
\]
which preserves the equation.

One possible approach is to derive a monotonicity formula for the energy of the extension and then adapt Quittner's argument~\cite{Quittner2021}. Here, we instead provide a more direct and simpler proof.

The proof first treats bounded entire solutions. The key step is to bound the local time dissipation of the weighted parabolic extension from below in terms of a positive value of the solution. Otherwise, after rescaling, compactness would produce a nontrivial stationary limit, contradicting the elliptic Liouville theorem. A covering argument then controls the measures of superlevel sets by the time dissipation. Integrating over the levels and combining the resulting estimate with a local energy inequality gives an iteration that forces the dissipation to vanish in the subcritical range. This contradicts the lower bound and proves nonexistence for bounded solutions.

To obtain the local a priori estimate, we use the doubling argument of Pol\'a\v cik, Quittner, and Souplet~\cite{PQS2007}, rescaling, and local regularity estimates. Since local convergence alone does not justify passage to the limit in the nonlocal equation, we apply the convergence theorem of Chen, Guo, Li,
and Ouyang~\cite{CGLO2026} and exclude the possible nonnegative constant defect using the integral representation. Failure of the estimate would then produce a bounded nontrivial entire solution, contradicting the first step. Finally, applying the a priori estimate on expanding cylinders proves Theorem \ref{main-thm-1} without a boundedness assumption. The required local compactness follows from the estimates below, which depend only on local norms and require the equation only in \(Q_4\), under global nonnegativity.

\begin{proposition}\label{prop:local-regularity}
Let $u\in\mathcal L_s(\R^n\times\R)\cap L^\infty(Q_4)$ be nonnegative
in $\R^n\times\R$, and suppose that $H^s u=f$ distributionally in
$Q_4$, with $f\in L^\infty(Q_4)$.
For every $0<\theta<\min\{2s,1\}$, $u$ has a representative satisfying
\begin{equation}\label{eq:local-holder}
 \|u\|_{C^{\theta,\theta/2}_{x,t}(Q_2)}
 \leq C\left(
 \|u\|_{L^\infty(Q_4)}+\|f\|_{L^\infty(Q_4)}
 \right),
\end{equation}
where $C=C(n,s,\theta)>0$.
If, in addition, $f\in C^{\eta,\eta/2}_{x,t}(Q_4)$ for some
$0<\eta<\min\{1,2-2s\}, 2s+\eta\neq1,$
then
\begin{equation}\label{eq:local-schauder}
 \|u\|_{C^{2s+\eta,\,(2s+\eta)/2}_{x,t}(Q_2)}
 \leq C\left(
 \|u\|_{L^\infty(Q_4)}
 +\|f\|_{C^{\eta,\eta/2}_{x,t}(Q_4)}
 \right),
\end{equation}
where $C=C(n,s,\eta)>0$.
\end{proposition}

The distributional equation in
Proposition~\ref{prop:local-regularity} means that
\[
 \int_{\R^n\times\R}u(x,t)(-\partial_t-\Delta)^s\psi(x,t)
 \dd x\dd t
 =\int_{Q_4}f(x,t)\psi(x,t)\dd x\dd t
 \qquad\text{for every }\psi\in C_c^\infty(Q_4),
\]
where
\[
 (-\partial_t-\Delta)^s\psi(x,t)
 =\frac{1}{|\Gamma(-s)|}
 \int_0^\infty
 \bigl[\psi(x,t)-e^{r\Delta}\psi(x,t+r)\bigr]
 \frac{\dd r}{r^{1+s}},
\]
and $e^{r\Delta}$ acts on the spatial variables.
The tail condition in \eqref{eq:tail-space} ensures that the pairing
with compactly supported test functions is well defined.

\begin{remark}
Chen, Guo, and Li~\cite[Theorems~1 and~2]{CGL2026} obtained
interior H\"older and Schauder estimates in terms of local norms
for nonnegative classical entire solutions of
$(\partial_t-\Delta)^s u=f$ with $f\ge0$.
Proposition~\ref{prop:local-regularity} retains this dependence
on local norms under the assumption that $u$ is globally
nonnegative, while the equation is required only in $Q_4$,
in the distributional sense, and $f$ may change sign there.
This local formulation applies directly in the proof of
Theorem~\ref{main-thm-2}, where the rescaled equations hold
on expanding cylinders.
\end{remark}

The paper is organized as follows. Section 2 presents the integral representation, the properties of the extension, and the elliptic nonexistence and compactness results used in the proof. Section 3 proves the Liouville theorem for bounded entire solutions. Section 4 establishes the local regularity estimates and the universal pointwise bound, and completes the proof of the Liouville theorem without a boundedness assumption.

\section{Preliminaries}\label{preliminary}

We begin by collecting the notation used below and specifying
the weak formulation of the extension problem.

\medskip
\noindent\textbf{Notation.}
\begin{itemize}
\item
We write $X=(x,y)$ and set
\[
 \mathbb R^{n+1}_+=\mathbb R^n\times(0,\infty),
 \qquad a=1-2s\in(-1,1),
 \qquad d_s=\frac{2^{1-2s}\Gamma(1-s)}{\Gamma(s)}.
\]
For extension functions, $\nabla=\nabla_{x,y}$ and
$\operatorname{div}$ act on $(x,y)$, while $\Delta=\Delta_x$.

\item
Let $B_R(x_0)$ denote the open Euclidean ball. For
$P=(x_0,t_0)$, define
\[
 \begin{aligned}
 Q_R(P)&=B_R(x_0)\times(t_0-R^2,t_0+R^2),\\
 Q_R^+(P)&=B_R(x_0)\times(0,R)
                  \times(t_0-R^2,t_0+R^2).
 \end{aligned}
\]
We omit the center when it is the origin and set
$D_R=B_R\times(0,R)$.
The flat boundary of $Q_R^+(P)$ is identified with $Q_R(P)$.
We use the parabolic distance
\[
 d_{\mathrm p}((x,t),(z,\tau))
 =\max\{|x-z|,|t-\tau|^{1/2}\},
 \qquad
 \operatorname{dist}_{\mathrm p}(P,E)
 =\inf_{Q\in E}d_{\mathrm p}(P,Q).
\]

\item
The Gaussian kernel is
\[
 G_r(z)=(4\pi r)^{-n/2}e^{-|z|^2/(4r)},
 \qquad z\in\mathbb R^n,\quad r>0.
\]

\item
The weighted Sobolev space $H^1(D_R,y^a)$ has norm
\[
 \|V\|_{H^1(D_R,y^a)}^2
 =\int_{D_R}\bigl(|V|^2+|\nabla V|^2\bigr)y^a\,dx\,dy,
\]
where the derivatives are understood weakly.
The trace $V(\cdot,0)$ is the
$L^2_{\mathrm{loc}}(B_R)$ limit as $y\downarrow0$;
it agrees with the ordinary boundary values when $V$
is continuous up to $y=0$.

\item
We use the usual parabolic H\"older spaces
$C^{\gamma,\gamma/2}$, with spatial variables $x$ on $Q_R$
and $(x,y)$ on $Q_R^+$.
For $0<\gamma<1$, the H\"older condition uses the parabolic
distance. For $1<\gamma<2$, the function is
$\gamma/2$-H\"older in time and its first spatial derivatives
belong to $C^{\gamma-1,(\gamma-1)/2}$.

\item
The space
$C_c^\infty(\overline{\mathbb R^{n+1}_+}\times\mathbb R)$
consists of restrictions of smooth compactly supported
functions on the full space. These tests may be nonzero
on $y=0$.
\end{itemize}
The notation $|E|$ denotes Lebesgue measure.
Constants $C$ may change from line to line and depend only
on $n,s,p$, unless further dependence is indicated.

For locally bounded $g$, a local weighted energy solution of
\[
 y^aV_t=\operatorname{div}(y^a\nabla V)
 \quad\text{in }Q_R^+,
 \qquad -y^aV_y=g\quad\text{on }y=0,
\]
is a function satisfying
\[
 V\in L^2\bigl(I;H^1(D_\rho,y^a)\bigr)
       \cap L^\infty\bigl(I;L^2(D_\rho,y^a)\bigr)
\]
for every $0<\rho<R$ and every open interval
$I\Subset(-R^2,R^2)$, together with
\[
 \int_{Q_R^+}
 \bigl(-V\phi_t+\nabla V\cdot\nabla\phi\bigr)y^a\,dx\,dy\,dt
 =
 \int_{Q_R}g(x,t)\phi(x,0,t)\,dx\,dt
\]
for every smooth test $\phi$ whose support in $y\ge0$ is
a compact subset of
$B_R\times[0,R)\times(-R^2,R^2)$.
On $\mathbb R^{n+1}_+\times\mathbb R$, these conditions
are required on every $Q_R^+$.
 
 We first recall the equivalent  integral representation of the master equation, as introduced in \cite{CGM2025} and \cite{CGL2026}.
\begin{lemma}[Integral representation]\label{lem:representation}
Let $u\geq0$ be bounded on $\R^n\times \R$ and belonging to  $C^{\gamma,\gamma/2}_{loc}(\R^n\times \R)$ for some $\gamma>2s$. If $(\partial_t-\Delta)^s u=f\geq0$ is bounded, then
\begin{equation}\label{eq:representation}
u=c+\frac1{\Gamma(s)}\int_0^\infty r^{s-1}S_rf\dd r,\qquad c\geq0.
\end{equation}
In particular, for $f=u^p$, one has $c=0$.
\end{lemma}

We will use the extension characterization of $(\partial_t-\Delta)^s$
from~\cite{NS2016,ST2017} to express the master equation as a local parabolic
equation with a nonlinear boundary condition.

\begin{lemma}\label{lem:extension}
Let $0<s<1$, $p>1$, and let $u$ be a bounded nonnegative classical
solution of
\[
 (\partial_t-\Delta)^s u=u^p
 \quad\text{in }\mathbb R^n\times\mathbb R.
\]
Set $a=1-2s$ and $B=\|u\|_\infty$. For $y>0$, define
\begin{equation}\label{eq:sourceextension}
 U(x,y,t)=\frac1{\Gamma(s)}
 \int_0^\infty\int_{\mathbb R^n}
 \frac{r^{s-1}e^{-(|x-z|^2+y^2)/(4r)}}{(4\pi r)^{n/2}}
 u(z,t-r)^p\,dz\,dr,
\end{equation}
and set $U(x,0,t)=u(x,t)$. Then the following properties hold.
\begin{enumerate}[label=\textup{(\roman*)},leftmargin=*]
\item
The integral is finite, $U$ is continuous up to $y=0$ and smooth for
$y>0$, and
\begin{equation}\label{eq:vertical}
 0\le U(x,y,t)\le u(x,t)\le B,\qquad
 0\le-y^aU_y(x,y,t)\le d_sB^p\quad(y>0).
\end{equation}
\item
The derivative $u_t$ exists, and both $u_t$ and $U_t$ are bounded and
continuous. The function $U_t$ extends continuously to $y=0$, with
\begin{equation}\label{eq:extension-time-boundary}
 U_t(x,0,t)=u_t(x,t),\qquad
 \lim_{y\downarrow0}
 \|U_t(\cdot,y,\cdot)-u_t\|_{L^\infty(\mathbb R^n\times\mathbb R)}=0.
\end{equation}
For every $R>0$ and every bounded open interval $I$,
\begin{equation}\label{eq:extension-strong-energy-space}
 \int_I\int_{D_R}
 \bigl(U^2+U_t^2+|\nabla U|^2+|\nabla U_t|^2\bigr)
 y^a\,dx\,dy\,dt<\infty.
\end{equation}
The maps $t\mapsto U(\cdot,t)$ and $t\mapsto\nabla U(\cdot,t)$
are absolutely continuous on $I$ with values in $L^2(D_R,y^a\,dx\,dy)$ and
$L^2(D_R,y^a\,dx\,dy)^{n+1}$, respectively; their time derivatives
are $U_t$ and $\nabla U_t$.
\item
The interior equation holds pointwise,
\begin{equation}\label{eq:extension}
 y^aU_t=\operatorname{div}(y^a\nabla U)\quad\text{for }y>0,
\end{equation}
and
\begin{equation}\label{eq:extension-boundary-limit}
 \lim_{y\downarrow0}
 \|-y^aU_y(\cdot,y,\cdot)-d_su^p\|_{L^\infty(\mathbb R^n\times\mathbb R)}=0.
\end{equation}
Consequently, for every
$\phi\in C_c^\infty(\overline{\mathbb R^{n+1}_+}\times\mathbb R)$,
\begin{equation}\label{eq:weak-extension}
 \int_{\mathbb R}\int_{\mathbb R^{n+1}_+}
 \bigl(-U\phi_t+\nabla U\cdot\nabla\phi\bigr)y^a\,dx\,dy\,dt
 =d_s\int_{\mathbb R}\int_{\mathbb R^n}
 u^p\phi(x,0,t)\,dx\,dt.
\end{equation}
\end{enumerate}
\end{lemma}

\begin{proof}
Put $f=u^p$. Lemma~\ref{lem:representation} gives
\[
 u(x,t)=\frac1{\Gamma(s)}\int_0^\infty\int_{\mathbb R^n}
 \frac{r^{s-1}e^{-|x-z|^2/(4r)}}{(4\pi r)^{n/2}}
 f(z,t-r)\,dz\,dr.
\]
Since $f\ge0$ and $0<e^{-y^2/(4r)}\le1$, the integral defining $U$
is finite and $0\le U\le u\le B$. Differentiation in $y>0$ gives
\begin{equation}\label{eq:extension-normal-integral}
 -y^aU_y(x,y,t)=\frac{y^{2-2s}}{2\Gamma(s)}
 \int_0^\infty\int_{\mathbb R^n}
 \frac{r^{s-2}e^{-(|x-z|^2+y^2)/(4r)}}{(4\pi r)^{n/2}}
 f(z,t-r)\,dz\,dr.
\end{equation}
The kernel is nonnegative. Integrating first in $z$ and then using
\[
 \int_0^\infty r^{s-2}e^{-y^2/(4r)}\,dr
 =\left(\frac{y^2}{4}\right)^{s-1}\Gamma(1-s)
\]
shows that its total mass is $d_s$. This proves
$0\le-y^aU_y\le d_sB^p$.
Monotone convergence gives $U(x,y,t)\to u(x,t)$ as $y\downarrow0$.
Integrating the derivative bound in $y$ then gives
\begin{equation}\label{eq:extension-trace-bound}
 0\le u(x,t)-U(x,y,t)\le\frac{d_s}{2s}B^py^{2s}.
\end{equation}

We also have the integral formula
\begin{equation}\label{eq:extension-poisson-integral}
 U(x,y,t)=\frac{y^{2s}}{4^s\Gamma(s)}
 \int_0^\infty\int_{\mathbb R^n}
 \frac{e^{-(|x-z|^2+y^2)/(4r)}}{r^{1+s}(4\pi r)^{n/2}}
 u(z,t-r)\,dz\,dr.
\end{equation}
Indeed, substitute the integral representation of $u$ into the
right-hand side. Tonelli's theorem, the convolution identity for
Gaussian kernels, and
\[
 \int_0^\ell e^{-y^2/(4r)}r^{-1-s}(\ell-r)^{s-1}\,dr
 =\Gamma(s)\left(\frac4{y^2}\right)^s
 \ell^{s-1}e^{-y^2/(4\ell)}
\]
give \eqref{eq:sourceextension}. 
The kernel in \eqref{eq:extension-poisson-integral} has mass one.
Writing $\tau=t-r$, all derivatives in $(x,y,t)$ fall on the kernel;
they are integrable locally uniformly for $y>0$, and the contribution
at $r=0$ vanishes because of $e^{-y^2/(4r)}$.
Differentiation therefore proves interior smoothness and
\eqref{eq:extension}. Together with \eqref{eq:extension-trace-bound},
this proves (i).

We next establish the time regularity in (ii). For bounded distributional
solutions of $(\partial_t-\Delta)^s w=F$ on all space-time, we use
\begin{equation}\label{eq:extension-holder-input}
 \sup_P\|w\|_{C^{\theta,\theta/2}(Q_1(P))}
 \le C_{n,s,\theta}\bigl(\|w\|_\infty+\|F\|_\infty\bigr),
 \qquad 0<\theta<\min\{2s,1\}.
\end{equation}
This follows from \cite[Theorem~1.20]{ST2017} by localization and
mollification. To check the operator-domain hypothesis there, choose
$\zeta\in C_c^\infty(Q_3(P))$ with $\zeta=1$ on $Q_{5/2}(P)$.
In $Q_2(P)$, the equation for $\zeta w$ has right-hand side $F+T$, where
\[
 T(x,t)=\frac1{|\Gamma(-s)|}
 \int_0^\infty\int_{\mathbb R^n}
 \frac{e^{-|x-z|^2/(4r)}(1-\zeta(z,t-r))w(z,t-r)}
 {r^{1+s}(4\pi r)^{n/2}}\,dz\,dr.
\]
Support separation for small $r$ and integrability for large $r$ give
$\|T\|_{L^\infty(Q_2(P))}\le C\|w\|_\infty$.
Mollifying $\zeta w$ gives smooth compactly supported functions in the
$L^2$ domain of $(\partial_t-\Delta)^s$, with uniform source bounds
on a slightly smaller cylinder. The cited estimate and passage to the
limit give \eqref{eq:extension-holder-input}.
When $2s=1$, the Zygmund estimate implies every exponent $\theta<1$.

Choose $0<\kappa<\min\{s,1/2\}$ with $1/\kappa\notin\mathbb N$, and put
$\delta_hu(x,t)=u(x,t+h)-u(x,t)$.
Estimate~\eqref{eq:extension-holder-input} with $\theta=2\kappa$
first gives $\|\delta_hu\|_\infty\le C|h|^\kappa$ for $|h|\le1$.
Suppose this bound holds with an exponent $\beta<1$.
Translation invariance gives
\[
 w_h=\frac{\delta_hu}{|h|^\beta},\qquad
 (\partial_t-\Delta)^s w_h=A_hw_h,\qquad
 A_h=p\int_0^1\bigl((1-\lambda)u(x,t)+\lambda u(x,t+h)\bigr)^{p-1}
 \,d\lambda.
\]
Since $0\le A_h\le pB^{p-1}$ and $\|w_h\|_\infty\le C$,
another application of \eqref{eq:extension-holder-input} gives
\[
 \|\delta_h^2u\|_\infty\le C|h|^{\beta+\kappa}.
\]
For $q=\beta+\kappa<1$, the identity
\[
 \delta_hu=2^{-N}\delta_{2^Nh}u
 -\sum_{k=0}^{N-1}2^{-k-1}\delta_{2^kh}^2u,
 \qquad \tfrac12\le2^N|h|\le1,
\]
implies $\|\delta_hu\|_\infty\le C|h|^q$ by summing a geometric series.
Finite iteration therefore gives
$\|\delta_h^2u\|_\infty\le C|h|^{1+\delta}$ for some $\delta>0$.
For $D_hu=\delta_hu/h$, this implies
\[
 \|D_hu-D_{h/2}u\|_\infty\le C|h|^\delta.
\]
Thus $D_{2^{-k}}u$ converges uniformly. It also converges distributionally
to $\partial_tu$, so its limit is the continuous derivative $u_t$.
The initial H\"older estimate makes each difference quotient uniformly
continuous in $(x,t)$; hence $u_t$ is bounded and uniformly continuous.
The same is true of $f_t=pu^{p-1}u_t$.

Differentiating \eqref{eq:extension-poisson-integral} and
\eqref{eq:extension-normal-integral} in time gives
\begin{align}
 U_t(x,y,t)
 &=\frac{y^{2s}}{4^s\Gamma(s)}
 \int_0^\infty\int_{\mathbb R^n}
 \frac{e^{-(|x-z|^2+y^2)/(4r)}}{r^{1+s}(4\pi r)^{n/2}}
 u_t(z,t-r)\,dz\,dr,
 \label{eq:extension-ut-integral}\\
 -y^a(U_t)_y(x,y,t)
 &=\frac{y^{2-2s}}{2\Gamma(s)}
 \int_0^\infty\int_{\mathbb R^n}
 \frac{r^{s-2}e^{-(|x-z|^2+y^2)/(4r)}}{(4\pi r)^{n/2}}
 f_t(z,t-r)\,dz\,dr.
 \label{eq:extension-ut-normal-integral}
\end{align}
These differentiations are justified by boundedness of $u_t,f_t$ and
integrability of the kernels. Their respective masses are $1$ and $d_s$,
so
\[
 \|U_t\|_\infty\le\|u_t\|_\infty,\qquad
 \|y^a(U_t)_y\|_\infty\le d_s\|f_t\|_\infty.
\]
In \eqref{eq:extension-ut-integral} and
\eqref{eq:extension-normal-integral}, set $r=y^2\tau$ and $z=x-y\xi$.
The resulting kernels in $(\xi,\tau)$ are independent of $y$ and
integrable. Uniform continuity of $u_t$ and $f$, first on a bounded
set of $(\xi,\tau)$ and then using the small integral of the kernel
outside that set, proves \eqref{eq:extension-time-boundary} and
\eqref{eq:extension-boundary-limit}.

Fix $R>0$ and a bounded open interval $I$. Choose a bounded open
interval $J$ with $\overline I\subset J$ and a smooth cutoff
$0\le\eta\le1$ such that
\[
 \eta=1\quad\text{on }\overline D_R\times\overline I,\qquad
 \operatorname{supp}(\eta|_{y\ge0})
 \subset B_{2R}\times[0,2R)\times J.
\]
For $0<\varepsilon<R$, test the classical equation for $U_t$ on
$y>\varepsilon$ with $\eta^2U_t$. Integration by parts gives
\[
 \begin{aligned}
 \int_{y>\varepsilon}\eta^2|\nabla U_t|^2y^a\,dx\,dy\,dt
 &=\int_{y>\varepsilon}
 \bigl(\eta\eta_tU_t^2-2\eta U_t\nabla U_t\cdot\nabla\eta\bigr)
 y^a\,dx\,dy\,dt\\
 &\quad+\int_{\mathbb R^n\times\mathbb R}
 \bigl[\eta^2U_t(-\varepsilon^a(U_t)_y)\bigr]_{y=\varepsilon}
 \,dx\,dt.
 \end{aligned}
\]
The absolute value of the boundary integral is at most
$C_\eta d_s\|u_t\|_\infty\|f_t\|_\infty$.
Young's inequality absorbs half the gradient integral. The remaining
terms are bounded independently of $\varepsilon$, since $a>-1$ and
$U_t$ is bounded. Letting $\varepsilon\downarrow0$ gives
$\nabla U_t\in L^2(D_R\times I,y^a\,dx\,dy\,dt)$.
The same calculation with $\eta^2U$, using $U\le B$ and
$0\le-y^aU_y\le d_sB^p$, gives the corresponding bound for
$\nabla U$. This proves \eqref{eq:extension-strong-energy-space}.
The distributional identities $\partial_tU=U_t$ and
$\partial_t\nabla U=\nabla U_t$, together with these $L^2$ bounds,
give the asserted absolute continuity in time. The continuous
representatives agree with the classical interior values by interior
smoothness. Integrating the equation for $U$ on $y>\varepsilon$ against a smooth
test and passing to the limit by the energy bounds and
\eqref{eq:extension-boundary-limit} proves \eqref{eq:weak-extension}.
\end{proof}

Note that when $u$ and $U$ are independent of $t$, then the extension problem of $U$ is reduced to a elliptic extension problem,  whose Liouville theorem is classical. 

\begin{lemma}[Elliptic obstruction]\label{lem:elliptic}
Let $n\ge2$ and $1<p<(n+2s)/(n-2s)$.
There is no bounded nonnegative function
$V\in H^1_{\mathrm{loc}}(\overline{\mathbb R^{n+1}_+},y^a)$
with a nonzero continuous trace $v$, satisfying
\[
 \operatorname{div}(y^a\nabla V)=0\quad\text{in }\mathbb R^{n+1}_+,
 \qquad -y^aV_y=d_sv^p\quad\text{on }y=0
\]
in the weak sense. Explicitly, the latter means
\[
 \int_{\mathbb R^{n+1}_+}y^a\nabla V\cdot\nabla\phi\,dxdy
 =d_s\int_{\mathbb R^n}v^p\phi(x,0)\,dx
 \quad\text{for every }\phi\in C_c^\infty(\overline{\mathbb R^{n+1}_+}).
\]
\end{lemma}

\begin{proof}
Let $Pv$ be the bounded Caffarelli--Silvestre Poisson extension of $v$:
\[
 Pv(x,y)=c_{n,s}\int_{\mathbb R^n}
 \frac{y^{2s}v(z)}{(|x-z|^2+y^2)^{(n+2s)/2}}\,dz,
\]
where the kernel has mass one. It is weighted harmonic and has trace
$v$ \cite{CS2007}. We first verify its local weighted energy.
The weighted trace theorem applied to localizations of $V$ gives
$v\in H^s_{\mathrm{loc}}(\mathbb R^n)$. Fix $R>0$ and choose
$\zeta\in C_c^\infty(\mathbb R^n)$ with $\zeta=1$ on $B_{2R}$.
Then $\zeta v\in H^s(\mathbb R^n)$, and the extension energy identity
\cite[Section~3.2]{CS2007} gives
\[
 \int_{\mathbb R^{n+1}_+}y^a|\nabla P(\zeta v)|^2\,dx\,dy
\le C_{n,s}\int_{\mathbb R^n}
|\xi|^{2s}|\widehat{\zeta v}(\xi)|^2\,d\xi<\infty.
\]
For $F=P((1-\zeta)v)$, separation of $B_R$ from the support of
$(1-\zeta)v$, followed by differentiation of the Poisson kernel, yields
\[
 |\nabla_xF(x,y)|\le C_Ry^{2s},\qquad
 |F_y(x,y)|\le C_Ry^{2s-1}
 \quad\text{on }B_R\times(0,R).
\]
Since $a=1-2s$, these estimates imply
\[
 \int_{B_R\times(0,R)}y^a|\nabla F|^2\,dx\,dy
 \le C_R\int_0^R\bigl(y^{1+2s}+y^{2s-1}\bigr)\,dy<\infty.
\]
Together with boundedness of $Pv$ and $a>-1$, this proves
$Pv\in H^1_{\mathrm{loc}}(\overline{\mathbb R^{n+1}_+},y^a)$.

Set $w=V-Pv$. It is bounded, weighted harmonic, and has zero trace.
Define its odd reflection by
$W(x,y)=\operatorname{sgn}(y)w(x,|y|)$ for $y\ne0$.
To justify $W\in H^1_{\mathrm{loc}}(\mathbb R^{n+1},|y|^a)$ across
$y=0$, take a smooth cutoff $\chi_\varepsilon(y)$ equal to zero for
$y\le\varepsilon$ and to one for $y\ge2\varepsilon$, with
$|\chi_\varepsilon'|\le C/\varepsilon$.
The zero trace and weighted Cauchy--Schwarz inequality give, for almost
every $x$,
\[
 |w(x,y)|^2\le\frac{y^{1-a}}{1-a}
 \int_0^y r^a|w_y(x,r)|^2\,dr.
\]
Consequently, for every compact $K\subset\mathbb R^n$,
\[
 \varepsilon^{-2}\int_{K\times(0,2\varepsilon)}y^a|w|^2\,dx\,dy
 \le C\int_{K\times(0,2\varepsilon)}y^a|w_y|^2\,dx\,dy
 \longrightarrow0.
\]
Thus $\chi_\varepsilon w\to w$ locally in weighted $H^1$.

Fix a ball $B\subset\mathbb R^n$ and $h>0$, and set
$D=B\times(0,h)$ and $\widetilde D=B\times(-h,h)$.
The preceding estimate gives
$w_\varepsilon:=\chi_\varepsilon w\to w$ in $H^1(D,y^a)$.
Let $W_\varepsilon$ be the odd reflection of $w_\varepsilon$.
Since $W_\varepsilon$ vanishes throughout $|y|<\varepsilon$, its two
halves agree with zero near $y=0$, so
$W_\varepsilon\in H^1(\widetilde D,|y|^a)$.
Reflection preserves the lengths of the gradients, and the weight is even.
Thus
\[
 \|W_\varepsilon-W_\delta\|_{H^1(\widetilde D,|y|^a)}^2
 =2\|w_\varepsilon-w_\delta\|_{H^1(D,y^a)}^2.
\]
It follows that $\{W_\varepsilon\}$ is Cauchy in $H^1(\widetilde D,|y|^a)$.
By completeness it has a limit there, and the identity
\[
 \|W_\varepsilon-W\|_{L^2(\widetilde D,|y|^a)}^2
 =2\|w_\varepsilon-w\|_{L^2(D,y^a)}^2\longrightarrow0
\]
identifies this limit with $W$. Since $B$ and $h$ were arbitrary,
$W\in H^1_{\mathrm{loc}}(\mathbb R^{n+1},|y|^a)$.

For $\phi\in C_c^\infty(\mathbb R^{n+1})$, put
$\psi(x,y)=\phi(x,y)-\phi(x,-y)$ for $y>0$.
Since $\psi(x,0)=0$, smoothness gives $|\psi(x,y)|\le C_\phi y$ near
$y=0$, and hence
\[
 \int_{\mathbb R^{n+1}_+}
 y^a|\nabla(\chi_\varepsilon\psi-\psi)|^2\,dx\,dy
 \le C_\phi\varepsilon^{a+1}\longrightarrow0.
\]
Testing the interior equation for $w$ with $\chi_\varepsilon\psi$
and passing to the limit is therefore justified by its local energy.
Changing variables in the lower half-space then gives
\[
 \int_{\mathbb R^{n+1}}|y|^a\nabla W\cdot\nabla\phi\,dx\,dy
 =\int_{\mathbb R^{n+1}_+}y^a\nabla w\cdot\nabla\psi\,dx\,dy=0.
\]
Thus $W$ is a bounded global weak energy solution of
$\operatorname{div}(|y|^a\nabla W)=0$.
Since $-1<a<1$, the weight $|y|^a$ belongs to $A_2$, and the coefficient
matrix $|y|^aI$ satisfies the weighted ellipticity condition with
constant one. The interior H\"older estimate
\cite[Theorem~2.3.12]{FKS} therefore applies and gives
\[
 \operatorname{osc}_{B_r}W\le C(r/R)^\delta\|W\|_\infty,
 \qquad 0<r<R/2,
\]
where the balls lie in $\mathbb R^{n+1}$ and $C,\delta>0$ depend only
on $n,s$. Fixing $r$ and letting $R\to\infty$ shows that $W$ is
constant. Oddness then gives $W=0$, hence $V=Pv$.

The distributional Dirichlet-to-Neumann identity \cite{CS2007}, compared
with the assumed weak boundary condition, now yields
\[
 \int_{\mathbb R^n}v(-\Delta)^s\varphi\,dx
 =\int_{\mathbb R^n}v^p\varphi\,dx
 \quad\text{for every }\varphi\in C_c^\infty(\mathbb R^n).
\]
Thus $(-\Delta)^sv=v^p\ge0$ distributionally. 
Since $v$ is continuous and nonnegative in all of $\mathbb R^n$, the
strong maximum principle \cite[Proposition~2.17]{Sil}, applied in a
ball containing any hypothetical zero of $v$, would imply $v\equiv0$.
This is excluded by the nonzero-trace assumption, so $v>0$ everywhere.
Finally, since $n>2s$, standard elliptic regularity and
Lemma~\ref{lem:representation}, applied to $v$ as a stationary function, give
\[
 v(x)=\kappa_{n,s}\int_{\mathbb R^n}
 \frac{v(z)^p}{|x-z|^{n-2s}}\,dz,
 \qquad \kappa_{n,s}>0.
\]
After a constant rescaling of $v$, the subcritical nonexistence
theorem of Chen--Li--Ou \cite{CLO} yields a contradiction.
\end{proof}

Therefore, if $u$ is a nontrivial nonnegative solution of $\eqref{master-lane-emden}$, it is natural to show that the local weighted energy $\int_{Q_R^+}y^a U^2_t \dd x\dd y \dd t$  is bounded away from zero. We will use the compactness argument to lead to the elliptic extension problem. Hence, the following boundary compactness of the extension problem is needed.

\begin{lemma}[Boundary compactness]\label{lem:compact}
Let $V_j$ be local weighted energy solutions of
\[
 y^a(V_j)_t=\operatorname{div}(y^a\nabla V_j)\quad\text{in }Q_2^+,
 \qquad -y^a(V_j)_y=g_j\quad\text{on }y=0,\ (x,t)\in Q_2,
\]
with the boundary condition understood in the weak sense. Assume that
\[
 \|V_j\|_{L^\infty(Q_2^+)}+\|g_j\|_{L^\infty(Q_2)}\le M
 \qquad\text{for all }j.
\]
Then each $V_j$ has a continuous representative up to $y=0$, and for every
$0<\rho<2$ and $0<\theta<\min\{2s,1\}$,
\[
 \|V_j\|_{C^{\theta,\theta/2}(\overline{Q_\rho^+})}
       \le C_{n,s,\rho,\theta}M,\qquad
 \int_{Q_\rho^+}|\nabla V_j|^2 y^a\,dx\,dy\,dt
       \le C_{n,s,\rho}M^2 .
\]
There exist a single subsequence $\{V_{j_k}\}$ and a function $V$ such
that, for every $0<\rho<2$,
\[
 \begin{aligned}
 V_{j_k}&\longrightarrow V
   &&\text{uniformly on }
      \overline{Q_\rho^+}=\overline B_\rho\times[0,\rho]\times[-\rho^2,\rho^2],\\
 \nabla V_{j_k}&\rightharpoonup\nabla V
   &&\text{in }L^2(Q_\rho^+,y^a\,dx\,dy\,dt)^{n+1}.
 \end{aligned}
\]
In particular, the traces $V_{j_k}(x,0,t)$ converge uniformly to
$V(x,0,t)$ on $\overline Q_\rho$. If $g_j\to g$ uniformly on every
compact subset of $Q_2$, then $V$ satisfies
\begin{equation}\label{eq:compact-limit-weak}
 \int_{Q_2^+}\bigl(-V\phi_t+\nabla V\cdot\nabla\phi\bigr)
                 y^a\,dx\,dy\,dt
   =\int_{Q_2}g(x,t)\phi(x,0,t)\,dx\,dt
\end{equation}
for every $\phi\in C_c^\infty(B_2\times(-2,2)\times(-4,4))$,
restricted to $y\ge0$.
\end{lemma}

\begin{proof}
We write $V_j=\Phi_j+Z_j$, choosing $\Phi_j$ to have boundary term
$\int_{Q_2}g_j\phi(x,0,t)\,dx\,dt$. Set $F_j=g_j/d_s$ on $Q_2$ and $F_j=0$ outside $Q_2$. Define
\[
 P_y(z,r)=\frac{r^{s-1}}{\Gamma(s)}G_r(z)e^{-y^2/(4r)},\qquad
 \Phi_j(x,y,t)=\int_0^\infty\!\int_{\mathbb R^n}
 P_y(z,r)F_j(x-z,t-r)\,dz\,dr .
\]
For $|t|\le4$, only $0<r<8$ contributes. Thus $\|G_r\|_1=1$ gives
$|\Phi_j(x,y,t)|\le CM$ for $y\ge0$ and $|t|\le4$.

For smooth $F_j$, the identity
$\partial_rP_y=(\Delta_z+\partial_{yy}+a y^{-1}\partial_y)P_y$
gives the interior equation for $\Phi_j$. Differentiating in $y$ gives
\[
 -y^a\partial_y\Phi_j(x,y,t)
 =\frac{y^{2-2s}}{2\Gamma(s)}
 \int_0^\infty\int_{\mathbb R^n}
 \frac{e^{-(|x-z|^2+y^2)/(4r)}}
      {r^{2-s}(4\pi r)^{n/2}}
 F_j(z,t-r)\,dz\,dr.
\]
The kernel is nonnegative and has total mass $d_s$.
After division by $d_s$, these kernels form an approximate
identity in $(x,t)$. Hence
$-y^a\partial_y\Phi_j\to d_sF_j$ locally uniformly as $y\downarrow0$.

For bounded measurable $F_j$, both differentiations remain valid for
$y>0$ by writing the integral in the source variables and differentiating
the kernel. The kernels in the preceding formula, divided by $d_s$,
have mass one and form an approximate identity in $(x,t)$. Since $F_j\in L^1\cap L^\infty$, it follows that
\[
 \|-y^a\partial_y\Phi_j(\cdot,y,\cdot)\|_\infty\le d_s\|F_j\|_\infty,
 \qquad
 -y^a\partial_y\Phi_j(\cdot,y,\cdot)\longrightarrow d_sF_j
 \quad\text{in }L^1(\mathbb R^n\times\mathbb R).
\]

Fix $0<\rho<2$ and set $R=(\rho+2)/2$. Choose a smooth cutoff $0\le\eta\le1$,
depending only on $\rho$, with
\[
 \eta=1\quad\text{on }\overline{Q_\rho^+},\qquad
 \operatorname{supp}(\eta|_{y\ge0})
       \subset B_R\times[0,R)\times(-R^2,R^2).
\]
For $0<\varepsilon<\rho$, test the interior equation on $y>\varepsilon$
with $\eta^2\Phi_j$. Integration by parts in time and Young's inequality
for the term $2\eta\Phi_j\nabla\eta\cdot\nabla\Phi_j$ give
\begin{equation}\label{eq:compact-caccioppoli}
 \begin{aligned}
 \int_{Q_R^+\cap\{y>\varepsilon\}}
   \eta^2|\nabla\Phi_j|^2\,y^a\,dx\,dy\,dt
 &\le C\int_{Q_R^+}\Phi_j^2
   (|\nabla\eta|^2+|\eta\eta_t|)\,y^a\,dx\,dy\,dt\\
 &\quad+C\int_{Q_R}
   \bigl|\eta^2\Phi_j\,\varepsilon^a\partial_y\Phi_j\bigr|_{y=\varepsilon}
   \,dx\,dt
 \le C_{n,s,\rho}M^2 .
 \end{aligned}
\end{equation}
Letting $\varepsilon\downarrow0$ proves the gradient bound for $\Phi_j$
on $Q_\rho^+$. Since $\rho<2$ is arbitrary, $\Phi_j$ has local weighted
energy up to $y=0$. Integrating the equation on $y>\varepsilon$ against
a smooth test $\phi$ and letting $\varepsilon\downarrow0$, the boundary
term converges to $\int_{Q_2}g_j\phi(x,0,t)\,dx\,dt$ by the $L^1$
convergence above. Thus $\Phi_j$ satisfies \eqref{eq:compact-limit-weak} with
$(\Phi_j,g_j)$ in place of $(V,g)$.

We next estimate the modulus of continuity of $\Phi_j$.
It is easy to check that the following $L^1$ bounds in $z$ hold uniformly for $y\ge0$.
\[
 \|P_y(\cdot,r)\|_1\le Cr^{s-1},\qquad
 \|\nabla_{z,y}P_y(\cdot,r)\|_1\le Cr^{s-3/2},\qquad
 \|\partial_rP_y(\cdot,r)\|_1\le Cr^{s-2}.
\]

Fix a time $t$ and two spatial points $X=(x,y)$ and $X'=(x',y')$
with $\delta=|X-X'|<1$ in $\overline{D_2}$. In the spatial convolution,
keep the source point fixed. At each $r$, the $L^1$ norm of the
difference of the two kernels is bounded by $Cr^{s-1}$ by the triangle
inequality and by $C\delta r^{s-3/2}$ by the mean value theorem.
Using the former for $r<\delta^2$ and the latter for $r\ge\delta^2$
gives, for $|t|\le4$,
\[
 \begin{aligned}
 |\Phi_j(X,t)-\Phi_j(X',t)|
 &\le CM\left(\int_0^{\delta^2}r^{s-1}\,dr
          +\delta\int_{\delta^2}^{8}r^{s-3/2}\,dr\right)\\
 &\le CM
 \begin{cases}
   \delta^{2s},&s<\tfrac12,\\
   \delta(1+|\log\delta|),&s=\tfrac12,\\
   \delta,&s>\tfrac12.
 \end{cases}
 \end{aligned}
\]
For a time increment $0<h<1$, keep the source time $\tau$ fixed.
The new interval $t<\tau<t+h$ contributes at most
$CM\int_0^h r^{s-1}\,dr$. On $\tau<t$ the data are the same in both
integrals, and the kernel difference is $P_y(\cdot,r+h)-P_y(\cdot,r)$
with $r=t-\tau$. For $r<h$ its $L^1$ norm is at most $Cr^{s-1}$.
For $r\ge h$, the estimate for $\partial_rP_y$ gives $Ch r^{s-2}$.
Hence
\[
 |\Phi_j(X,t+h)-\Phi_j(X,t)|
 \le CM\left(\int_0^h r^{s-1}\,dr
              +h\int_h^\infty r^{s-2}\,dr\right)
 \le CMh^s.
\]
Combining the two estimates with the supremum bound yields
$\|\Phi_j\|_{C^{\theta,\theta/2}(\overline{Q_2^+})}\le C_\theta M$
for every $0<\theta<\min\{2s,1\}$. At $s=\tfrac12$ we use
$\delta(1+|\log\delta|)\le C_\theta\delta^\theta$ for $\theta<1$.

Set $Z_j=V_j-\Phi_j$. Subtracting their weak identities gives
\[
 \int_{Q_2^+}\bigl(-Z_j\partial_t\phi+\nabla Z_j\cdot\nabla\phi\bigr)
 y^a\,dx\,dy\,dt=0
\]
for every test $\phi$ in \eqref{eq:compact-limit-weak}.Testing the weak equation with $\eta^2 Z_j$ and integrating by parts
in time, we obtain
\[
 \int_{Q_R^+}\eta^2|\nabla Z_j|^2\,y^a\,dx\,dy\,dt
 =
 \int_{Q_R^+}
 \bigl(\eta\eta_t Z_j^2
       -2\eta Z_j\,\nabla Z_j\cdot\nabla\eta\bigr)
 y^a\,dx\,dy\,dt .
\]
The integration by parts follows from the energy-space chain rule
applied to $\eta Z_j$, which vanishes at both time endpoints.
By Young's inequality, the mixed term is bounded by half the
left-hand side plus
$2\int_{Q_R^+}Z_j^2|\nabla\eta|^2\,y^a\,dx\,dy\,dt$.
Since $\eta=1$ on $Q_\rho^+$ and $|Z_j|\le CM$, it follows that
\[
 \int_{Q_\rho^+}|\nabla Z_j|^2\,y^a\,dx\,dy\,dt
 \le C\int_{Q_R^+} Z_j^2
       \bigl(|\nabla\eta|^2+|\eta\eta_t|\bigr)
       y^a\,dx\,dy\,dt
 \le C_{n,s,\rho}M^2 .
\] Together with \eqref{eq:compact-caccioppoli},
this gives the asserted gradient bound for $V_j=\Phi_j+Z_j$.

Apply the boundary estimate of
Audrito--Fioravanti--Vita \cite[Theorem~1.1(i)]{AFV2024} on half-cylinders
centered at $y=0$. The weight satisfies $a>-1$, the coefficient matrix
is the identity, and both interior forcing terms vanish. Their integrability
exponents may therefore be chosen large enough for any fixed $\theta<1$.
The theorem bounds the $C^{\theta,\theta/2}$ norm on a smaller half-cylinder
by the weighted $L^2$ norm of $Z_j$ on the larger one. This last norm is
bounded by $C_{n,s,\rho}M$. Since $\rho<R$, a finite covering along $y=0$
and interior estimates away from $y=0$ give the bound on
$\overline{Q_\rho^+}$. Adding the estimate for $\Phi_j$ proves the claimed
H\"older bound for $V_j$.

Take $\rho_m=2-1/m$, $m\ge1$. On each $\overline{Q_{\rho_m}^+}$,
Arzel\`a--Ascoli gives a uniformly convergent subsequence, and the energy
bound gives a further subsequence whose gradients converge weakly in
$L^2(Q_{\rho_m}^+,y^a\,dx\,dy\,dt)^{n+1}$. Testing against smooth
functions supported in $y>0$ identifies this weak limit with $\nabla V$,
where $V$ is the uniform limit. Since $\{y=0\}$ has zero weighted
volume, this identifies the gradient almost everywhere on the whole
half-cylinder. A diagonal choice gives the same subsequence for all $m$,
and hence for every $\rho<2$. Finally, every test in
\eqref{eq:compact-limit-weak}, restricted to $y\ge0$, has support in
$\overline{Q_{\rho_m}^+}$ for some $m$. In the weak identity for $V_{j_k}$,
the time term converges by uniform convergence, the gradient term by
weak convergence, and the boundary integral by local uniform convergence
of $g_{j_k}$ to $g$. This proves \eqref{eq:compact-limit-weak}.
\end{proof}

\section{Proof of Theorem \ref{main-thm-1} for bounded solution}
In this section, we always assume that $u$ is a nonnegative classical solution of \eqref{master-lane-emden} and $U$ is its extension. Set 
\begin{align*}
A(R)&=\int_{Q_R}u^{p+1}\dd x\dd t,\qquad\  D(R)=\int_{Q_R^+}y^a U^2_t \dd x\dd y \dd t\\ M(R)&=\int_{Q^+_R}y^aU^2 \dd x \dd y \dd t,\quad  G(R)=\int_{Q^+_R}y^a|\nabla U|^2\dd x\dd y\dd t.
\end{align*}
Set
$$
\alpha=\frac{2s}{p-1},\quad \mu=2\alpha+2s-n.
$$
Note that $p<p(n,s)$ is equivalent to $\mu>0$. 
Assume \eqref{master-lane-emden} admits positive bounded solution, we will obtain contradiction by estimating the decay of $D(R)$ and showing the local positive lower bound of $D(R)$.

First, by testing the extension equation by suitable functions, we can control $D(R)$ by $A(4R)$.

\begin{lemma}[Local dissipation estimate]\label{lem:local-dissipation}
For every $R>0$,
\begin{equation}\label{eq:local-dissipation}
 D(R)\le CR^{-2}A(4R)+CR^{-\mu}.
\end{equation}
\end{lemma}

\begin{proof}
Fix $R>0$. We take cutoff functions $\chi=\chi(x,y)$ and $\zeta=\zeta(t)$ to satisfy 
$$
\chi\in C^\infty_c(B_{2R}\times[0,2R)),\quad \chi=1\quad\text{in }B_R\times [0,R],\quad 0\leq \chi\leq 1,\quad |\nabla\chi|\leq \frac{C}{R}.
$$
and
$$
\zeta\in C^\infty_c(-4R^2,4R^2),\quad\zeta=1\quad\text{on }[-R^2,R^2],\quad 0\leq \zeta\leq 1,\quad |\zeta'|\leq\frac{C}{R^2}.
$$
By Lemma \ref{lem:extension}, we can define 
$$
E_\chi(t)=\frac12\int_{\R^{n+1}_+}y^a\chi^2|\nabla U|^2dxdy
 -\frac{d_s}{p+1}\int_{\R^n}\chi_0^2u^{p+1}dx,
$$
where $\chi_0(\cdot)=\chi(\cdot,0)$. By Lemma~\ref{lem:extension}, $E_\chi$ is locally absolutely
continuous, and the following derivative identity holds
for almost every $t$.
Multiplying both sides of \eqref{eq:extension} by  $\chi^2U_t$  and integrating about $x,y$  over $\mathbb R^{n+1}_+$, we have
$$
\int_{\mathbb R^{n+1}_+}y^a\chi^2U^2_t\dd x\dd y= d_s \int_{\mathbb R^n}\chi_0^2u^p u_t dx-\int_{\mathbb R^{n+1}_+}y^a\chi^2\nabla U\cdot\nabla U_t\dd x\dd y-2\int_{\mathbb R^{n+1}_+}y^a\chi U_t\nabla U\cdot\nabla\chi dxdy.
$$
Hence,
\begin{align*}
E'_\chi(t)&=\int_{\R^{n+1}_+}y^a\chi^2\nabla U\cdot\nabla U_t
 -d_s\int_{\R^n}\chi_0^2u^{p}u_t dx\\
 &=-\int_{\mathbb R^{n+1}_+}y^a\chi^2U^2_t\dd x\dd y-2\int_{\mathbb R^{n+1}_+}y^a\chi U_t\nabla U\cdot\nabla\chi dxdy.
\end{align*}
Next, multiplying this equation by $\zeta^2$  and integrating about $t$  over $\mathbb R$, we get
\begin{align*}
\int y^a\zeta^2\chi^2U^2_t&=-\int E'_{\chi}\zeta^2-2\int y^a\zeta^2 \chi U_t\nabla U\cdot\nabla\chi\\
&=2\int \zeta\zeta' E_\chi-2\int (y^{\frac a2}\zeta \chi U_t)(y^{\frac a2}\zeta\nabla U\cdot\nabla\chi)
\end{align*}
Young’s inequality gives
$$
\int y^a\zeta^2\chi^2U^2_t\leq CR^{-2}(G(2R)+A(2R))+\frac12\int y^a\zeta^2\chi^2U^2_t.
$$
Thus, 
\begin{equation}\label{D-G-A-estimate}
D(R)\leq \int y^a\zeta^2\chi^2U^2_t\leq CR^{-2}(G(2R)+A(2R)).
\end{equation}

To estimate $G(2R)$, we test \eqref{eq:extension} by $U\chi^2\zeta^2$, then
\begin{align}\label{test-2-ineq}
\int y^a U U_t\chi^2\zeta^2=d_s\int u^{p+1} \chi_0^2\zeta^2-\int y^a|\nabla U|^2\chi^2\zeta^2-2\int y^a \chi \zeta^2U\nabla U\cdot\nabla\chi
\end{align}
Note that
$$
\int y^a U U_t\chi^2\zeta^2=\frac12\int y^a (U^2)_t\chi^2\zeta^2=-\int y^aU^2\chi^2\zeta\zeta'
$$
Using the Young inequality in \eqref{test-2-ineq} implies
\begin{align*}
\int y^a|\nabla U|^2\chi^2\zeta^2&=\int y^a U^2 \chi^2\zeta\zeta'+d_s\int u^{p+1} \chi_0^2\zeta^2-2\int (y^{\frac a2} U\zeta\nabla\chi) \cdot (y^{\frac a2}\chi\zeta\nabla U)\\
&\leq CR^{-2}M(2R)+CA(2R)+\frac12 \int y^a|\nabla U|^2\chi^2\zeta^2 
\end{align*}
Replacing $R$ by $2R$, we have
\begin{equation}\label{G-estimate}
G(2R)\leq CR^{-2} M(4R)+CA(4R).
\end{equation}

It remains to estimate $M(4R)$. Note that $2s+2\alpha=\alpha(p+1)$. By the elementary inequality $z^2\leq 1+z^{p+1}$ for $p>1$ and any $z\geq 0$, we have 
$$
R^{-2s}u^2=R^{-\alpha(p+1)}(R^\alpha u)^2\leq u^{p+1}+R^{-\alpha(p+1)}.
$$
Recall $U(x,y,t)\leq u(x,t)$ for any $y>0$, hence  
\begin{equation}\label{M-estimate}
\begin{aligned}
R^{-2}M(4R)&\leq R^{-2}\int y^au^2 =\frac{4^{2-2s}}{2-2s}\int_{Q_{4R}} R^{-2s}u^2\leq C(A(4R)+|Q_{4R}|R^{-\alpha(p+1)})\\
&\leq C(A(4R)+R^{2-\mu}) 
\end{aligned}
\end{equation}

Combining \eqref{D-G-A-estimate}, \eqref{G-estimate}, and \eqref{M-estimate}, we get \eqref{eq:local-dissipation}.
\end{proof}

Next, the following key estimate controls $u$ pointwise by the local weighted energy, which will  not only be used to obtain the estimate of the level set of $u$, but also give a local positive lower bound of $D(R)$. The main idea is based on the compactness argument and reduces the desired estimate to the Liouville Theorem for the elliptic extension problem.

\begin{lemma}\label{lem:value-cost}
There exist constants  $L\ge1$ and $c>0$ such that, for every $P$ with $u(P)>0$,
\begin{equation}\label{eq:value-cost}
 \int_{Q_{L u(P)^{-1/\alpha}}^+(P)}y^aU_t^2
 \ge c\,u(P)^{\mu/\alpha}.
\end{equation}
\end{lemma}
\begin{proof}
Step 1. We claim that there exist $L_0\geq 1$ and $c_0>0$ such that for every bounded solution $v$ of \eqref{master-lane-emden} satisfying $v(0,0)=1$ and $0\leq v\leq 2^\alpha$ in $Q_{2L_0}$, its extension $V$ holds
\begin{equation}\label{normal-estimate}
\int_{Q^+_{L_0}}y^aV^2_t\geq c_0. 
\end{equation}

If the claim is false, then for every integer $j\ge1$, there exist a solution $v_j$ with extension $V_j$ such that
\[
 v_j(0,0)=1,\qquad 0\le v_j\le2^\alpha\quad \hbox{in }Q_{2j},
 \qquad \int_{Q_j^+}y^a(V_j)_t^2<j^{-1}.
\]

By \eqref{eq:vertical}, $0\le V_j(x,y,t)\le v_j(x,t)\le 2^\alpha$
whenever $(x,t)\in Q_{2j}$ and $y>0$.
On $y=0$, the boundary data satisfy
$0\le g_j=d_s v_j^p\le d_s2^{\alpha p}$ in $Q_{2j}$.


Lemma~\ref{lem:compact} with a diagonal argument gives a
bounded nonnegative limit $v$ in $\R^n\times\R$ and $V$ in
$\R^{n+1}_+\times\R$ satisfying $v(0,0)=1$ and $0\le v\le2^\alpha$. In addition, for any $\phi\in C^\infty_c(\mathbb R^n\times(0,\infty)\times\mathbb R)$,
\[
 \left|\int y^aV_j\phi_t\right|
 =\left|\int y^a(V_j)_t\phi\right|
 \le\left(\int_{\operatorname{supp}\phi}y^a(V_j)_t^2\right)^{1/2}
     \left(\int y^a\phi^2\right)^{1/2}\longrightarrow0.
\]
Thus, $V_t=0$ for $y>0$, and continuity makes $v_t=0$ as well. In the limit equation satisfied by $V$, we take
$\phi(x,y,t)=\psi(x,y)\zeta(t)$ as a test function, where $\psi\in C^\infty_c(\overline{\R^{n+1}_+})$, $\zeta\in C_c^\infty(\R)$,
and $\int_\R\zeta=1$.
Note that $\int_\R\zeta'=0$, which implies
\[
 \int_{\R^{n+1}_+}y^a\nabla V\cdot\nabla\psi
 =d_s\int_{\R^n}v^p\psi(x,0).
\]
This is in contradiction with  Lemma~\ref{lem:elliptic} since  $v(0,0)=1$. Hence the claim holds.

Step 2. We use the doubling argument to obtain the desired estimate.

We set $m=u^{1/\alpha}$ and $k=2L_0$. A doubling argument will show that there exists a
point $Z$ such that
\begin{equation}\label{eq:doubling-core}
 m(Z)\ge m(P),\qquad d_{\mathrm p}(Z,P)\le\frac{2k}{m(P)},\qquad
 m\le2m(Z)\quad\hbox{in }Q_{k/m(Z)}(Z).
\end{equation}
Indeed, we start with $Z_0=P$.  If the last inequality fails at $Z_i$, then we can choose
$Z_{i+1}\in Q_{k/m(Z_i)}(Z_i)$ with $m(Z_{i+1})>2m(Z_i)$. Since $u$
is bounded, this procedure will stop in finitely many steps. Moreover, assume that we stop with the last inequality holds for $Z_j$, then 
$$
d_{\mathrm p}(Z_j,P)\leq \sum_{i=0}^{j-1}d_{\mathrm p}(Z_i,Z_{i+1})\leq \sum^{j-1}_{i=0}\frac{k}{m(Z_i)}\leq \sum^{j-1}_{i=0}\frac{k}{2^im(P)}\leq \frac{2k}{m(P)}.
$$
Hence, we prove \eqref{eq:doubling-core}.

Set $r=m(Z)^{-1}$ with the above $Z$. We write $Z=(z,\tau)$ and define
\[
 v(x,t)=r^\alpha u(z+rx,\tau+r^2t),\qquad
 V(x,y,t)=r^\alpha U(z+rx,ry,\tau+r^2t).
\]
Note that $\alpha(p-1)=2s$, therefore $V$ is the extension of $v$ and $V$ satisfies the same extension equation and the
boundary condition. A change of variables in
\eqref{eq:sourceextension} shows that $V$ is the canonical extension
of $v$. In addition, $v(0,0)=1$ and $v\le2^\alpha$ in $Q_{2L_0}$ by \eqref{eq:doubling-core}. By applying the claim in Step 1, we have
$$
c_0\leq \int_{Q^+_{L_0}}y^aV_t^2=r^{2\alpha+4-a-n-3}\int_{Q^+_{L_0r}(Z)}y^aU_t^2=r^\mu \int_{Q^+_{L_0r}(Z)}y^aU_t^2,
$$
Note that $\mu>0$, \eqref{eq:doubling-core} yields 
$$
\int_{Q^+_{L_0r}}y^aU_t^2\geq c_0r^{-\mu}=c_0m(Z)^{\mu}\geq c_0 m(P)^\mu =c_0 u(P)^{\mu/\alpha}.
$$
We also have $d_{\mathrm p}(Z,P)\leq 4L_0/m(P)$ and $L_0r\leq L_0/m(P)$ and hence $Q^+_{L_0r}(Z)\subset Q^+_{5L_0/m(P)}(P)$. Therefore, we obtain \eqref{eq:value-cost} with $L=5L_0$ and $c=c_0$.
\end{proof}

Now, we are ready to estimate $A(R)$ by $D(2R)$. The key is to control the measure of the level sets of $u$ based on a parabolic grid decomposition.

\begin{lemma}\label{lem:distribution}
Set
\[
 q=\frac{n+2+\mu}{\alpha}=p+1+\frac2\alpha.
\]
For every $R>0$ and every $h\ge L^\alpha R^{-\alpha}$,
\begin{equation}\label{eq:distribution}
 \bigl|\{u>h\}\cap Q_R\bigr|\le C D(2R)h^{-q}.
\end{equation}
In addition, for the same $h$,
\begin{equation}\label{eq:layer-integral}
 A(R)\le CR^{n+2}h^{p+1}+CD(2R)h^{-2/\alpha}.
\end{equation}
\end{lemma}

\begin{proof}
Fix $R>0$ and  $\ell=h^{-1/\alpha}$ with $h\ge L^\alpha R^{-\alpha}$.  For $\mathbf k=(k_1,\ldots,k_n)\in\mathbb Z^n$
and $m\in\mathbb Z$, define
\[
 \mathcal C_{\mathbf k,m}
 =\prod_{\nu=1}^n[k_\nu\ell,(k_\nu+1)\ell)
   \times[m\ell^2,(m+1)\ell^2).
\]
These cells are pairwise disjoint and satisfy
\[
 \mathbb R^n\times\mathbb R
 =\bigcup_{(\mathbf k,m)\in\mathbb Z^n\times\mathbb Z}
   \mathcal C_{\mathbf k,m},
 \qquad
 |\mathcal C_{\mathbf k,m}|=\ell^{n+2}.
\]
There are finitely many cells 
that meet $Q_R$ and hence we denote 
$$
\{\mathcal C_{\mathbf k,m}:\mathcal C_{\mathbf k,m}\cap\{u>h\}\cap Q_R\neq\emptyset\}=\{\mathcal C_i\}^{J}_{i=1}.
$$
For each $1\leq i\leq J$, we choose one point $P_i\in \mathcal C_{i}\cap\{u>h\}\cap Q_R$. Then, by 
Lemma~\ref{lem:value-cost}, we have
\begin{equation}\label{Pi-estimate}
 \int_{Q_{L u(P_i)^{-1/\alpha}}^+(P_i)}y^aU_t^2
 \ge c\,u(P_i)^{\mu/\alpha}
 \ge c\,h^{\mu/\alpha}.
\end{equation}
These cylinders $\{Q_{L u(P_i)^{-1/\alpha}}^+(P_i)\}^J_{i=1}$ have radii at most $L\ell\le R$ and hence are contained in
$Q_{2R}^+$. For any point $(x,y,t)\in Q^+_{2R}$, if it is in $Q_{L u(P_i)^{-1/\alpha}}^+(P_i)$, then
$P_i=(x_i,t_i)$ satisfies
\[
 |x_i-x|<L\ell,\qquad |t_i-t|<L^2\ell^2.
\]
This rectangle meets with at most
$C_*=(\lceil2L\rceil+3)^n(\lceil2L^2\rceil+3)$ grid cells.
Each cell supplies at most one center, so the overlap of the cylinders at $(x,y,t)$ is at most $C_*$. Summing \eqref{Pi-estimate} about $i$ from $1$ to $J$  gives
\[
 cJ h^{\mu/\alpha}\le C_*D(2R).
\]
Since $\{\mathcal C_i\}^{J}_{i=1}$ cover $\{u>h\}\cap Q_R$, we get
\[
 \bigl|\{u>h\}\cap Q_R\bigr|
 \le J\ell^{n+2}
 \le CD(2R)h^{-(n+2+\mu)/\alpha},
\]
which is \eqref{eq:distribution}.

Moreover, for $h\ge L^\alpha R^{-\alpha}$, \eqref{eq:distribution} yields
\begin{align*}
 A(R)
 &=(p+1)\int_0^\infty s^p\bigl|\{u>s\}\cap Q_R\bigr|\,ds\\
 &\le (p+1)|Q_R|\int_{0}^hs^{p}
ds   +CD(2R)\int_h^\infty s^{p-q}\,ds\\
 &\le CR^{n+2}h^{p+1}+CD(2R)h^{-2/\alpha},
\end{align*}
where the last inequality uses $q=p+1+2/\alpha$.
Hence, we obtain \eqref{eq:layer-integral}.

\end{proof}

We can now prove Theorem \ref{main-thm-1} for bounded solution.

\begin{proposition}\label{prop:bounded-zero}
Every bounded nonnegative entire solution in the subcritical range is
identically zero.
\end{proposition}

\begin{proof}
We apply \eqref{eq:layer-integral} in $Q_{4R}$
with $h=KR^{-\alpha}$, where $K\ge (L/4)^\alpha$ and therefore
$h\ge L^\alpha(4R)^{-\alpha}$. Note that $n+2-\alpha(p+1)=2-\mu$, we obtain
\[
 A(4R)\leq CR^{n+2}(KR^{-\alpha})^{p+1}+CD(8R)(KR^{-\alpha})^{-2/\alpha}= CK^{p+1}R^{2-\mu}
       +CK^{-2/\alpha}R^2D(8R).
\]
Substitute it into \eqref{eq:local-dissipation} gives
\begin{equation}\label{eq:absorption-start}
 D(R)\le CK^{-2/\alpha}D(8R)
       +C(K^{p+1}+1)R^{-\mu}.
\end{equation}
where the constant $C$ is independent of $K$ and $R$.
Fix $\delta=8^{-n}/2$ and then choose $K$ larger so that
$CK^{-2/\alpha}\le\delta$. Thus
\begin{equation}\label{eq:radius-iteration}
 D(R)\le\delta D(8R)+C_\delta R^{-\mu}.
\end{equation}
Iteration gives
\begin{equation}\label{iteration-D}
 D(R)\le\delta^jD(8^jR)
    +C_\delta R^{-\mu}\sum_{i=0}^{j-1}(\delta8^{-\mu})^i.
\end{equation}
Set $B=\|u\|_{L^\infty(\mathbb R^n\times \mathbb R)}$.
For the first term on the right, we use the boundedness of $u$ and \eqref{eq:local-dissipation} 
and hence
\[
 \delta^jD(8^jR)\leq C\delta^j [B^{p+1} (8^jR)^n+(8^jR)^{-\mu}]
 \le C B^{p+1}R^n(\delta8^n)^j
     +CR^{-\mu}(\delta8^{-\mu})^j\rightarrow0,
\]
as $j\rightarrow\infty$.
Here $\delta8^n=1/2$ and $\delta8^{-\mu}<1$.
Letting $j\to\infty$ in \eqref{iteration-D} yields
\begin{equation}\label{eq:dissipation-decay}
 D(R)\le\frac{C_\delta}{1-\delta8^{-\mu}}R^{-\mu}
       \le CR^{-\mu}.
\end{equation}

Assume $u(P)>0$ at some point $P$. 
Then for $R$ sufficiently large, 
\eqref{eq:value-cost} and \eqref{eq:dissipation-decay} give
\[
 0<c\,u(P)^{\mu/\alpha}
 \le\int_{Q_{L u(P)^{-1/\alpha}}^+(P)}y^aU_t^2
 \le D(R)\le CR^{-\mu}\longrightarrow0.
\]
This contradiction proves $u\equiv0$.
\end{proof}

\section{Proof of the a priori local estimate}\label{sec:universal-estimate}

We first prove Proposition~\ref{prop:local-regularity}, which gives
the compactness needed for the blow-up argument. We choose a smooth
cutoff $\chi$ supported in $Q_4$ and equal to one on $Q_3$, and write
$u=\chi u+(1-\chi)u$. The first term agrees with $u$ in $Q_3$,
while the second accounts for the contribution from outside $Q_3$.
A similar cutoff argument appears in the proof of
\cite[Proposition~2.1]{SWWY}. Since $u\ge0$ globally, testing the
localized equation at later times allows us to bound this exterior
contribution using only the norms of $u$ and $f$ on $Q_4$, even when
$f$ changes sign. After mollification, the interior estimates of
Stinga and Torrea~\cite[Theorem~1.20]{ST2017} give the stated
H\"older and Schauder bounds.

\begin{proof}[Proof of Proposition~\ref{prop:local-regularity}]
Set $L=\|u\|_{L^\infty(Q_4)}+\|f\|_{L^\infty(Q_4)}$ and assume
$L<\infty$. Choose $\chi\in C_c^\infty(Q_4)$ with $0\le\chi\le1$
and $\chi=1$ on $Q_3$, and write
\[
 w=\chi u,\qquad h=(1-\chi)u\ge0.
\]
Then $w$ is bounded and compactly supported, and $h=0$ in $Q_3$. Define
\[
 K(z,r)=\frac{r^{-\nu}e^{-|z|^2/(4r)}}
 {(4\pi)^{n/2}|\Gamma(-s)|},\qquad
 \nu=\frac n2+1+s,\qquad r>0,
\]
and
\[
 T(x,t)=\int_{-\infty}^t\!\int_{\R^n}
 K(x-z,t-\tau)h(z,\tau)\,dz\,d\tau\ge0.
\]
We first verify that $T$ is locally bounded in $Q_3$.
Fix a compact set $D\Subset Q_3$ and choose $0<\varepsilon<1$
such that
\[
 D\subset B_{3-\varepsilon}
       \times(-9+\varepsilon,9-\varepsilon).
\]
For $(x,t)\in D$ and $\tau<t$, put
$r=t-\tau$ and $\delta=10-t$. Then $1<\delta<19$, and
completing the square gives
\[
 -\frac{|z-x|^2}{4r}+\frac{|z|^2}{4(r+\delta)}
 =
 -\frac{\delta}{4r(r+\delta)}
  \left|z-x-\frac r\delta x\right|^2
 +\frac{|x|^2}{4\delta}.
\]
For $r\ge\varepsilon/6$, this identity and
$1+(r+\delta)^\nu\le C_\varepsilon r^\nu$ imply
\[
 K(x-z,r)
 \le C_\varepsilon
 \frac{e^{-|z|^2/[4(r+\delta)]}}
      {1+(r+\delta)^\nu}.
\]
If $0<r<\varepsilon/6$ and $h(z,\tau)\ne0$, then
$\tau\in(-9,9)$ and hence $|z|\ge3$, because $h=0$ in $Q_3$.
Thus $|z-x|\ge\varepsilon$ and
\[
 \left|z-x-\frac r\delta x\right|
 \ge |z-x|-\frac{r|x|}{\delta}
 \ge\frac{\varepsilon}{2}.
\]
The same comparison therefore holds, since the ratio of its
left-hand side to the kernel on its right-hand side is bounded by
\[
 C r^{-\nu}e^{-\varepsilon^2/(32r)}
 \le C_\varepsilon.
\]
Since $r+\delta=10-\tau$ and $0\le h\le u$, we conclude that
\[
 \sup_{(x,t)\in D}T(x,t)
 \le C_D
 \int_{-\infty}^{10}\int_{\mathbb R^n}
 \frac{u(z,\tau)e^{-|z|^2/[4(10-\tau)]}}
      {1+(10-\tau)^\nu}\,dz\,d\tau
 <\infty,
\]
where finiteness follows from $u\in\mathcal L_s$. 

We next justify the localized equation in distributions. For
$\psi\in C_c^\infty(Q_3)$, let $H^*=-\partial_t-\Delta_x$.
Since $h\psi=0$, the adjoint semigroup formula gives, on the support
of $h$,
\[
 (H^*)^s\psi(z,\tau)
 =-\int_\tau^\infty\!\int_{\R^n}
 K(x-z,t-\tau)\psi(x,t)\,dx\,dt.
\]
Fubini is justified by $\int T|\psi|<\infty$, and hence
\[
 \langle H^s h,\psi\rangle
 =\int h(H^*)^s\psi=-\int T\psi.
\]
Together with the assumed equation for $u$, this yields
\begin{equation}\label{eq:cutoff-tail}
 H^s w=f+T\quad\text{in }\mathcal D'(Q_3),
 \qquad
 \int w(H^*)^s\psi=\int(f+T)\psi.
\end{equation}
All pairings on the left are over $\R^{n+1}$; the test functions
on the right are supported in $Q_3$.

Choose a fixed $\phi\in C_c^\infty(A)$ with $\phi\ge0$ and
$\int\phi=1$, where $A=B_1\times(7,8)\Subset Q_3$.
Testing \eqref{eq:cutoff-tail} gives
\begin{equation}\label{eq:tail-average}
 0\le\int\phi T
 =\int w(H^*)^s\phi-\int f\phi\le CL.
\end{equation}
Here $\|w\|_\infty\le\|u\|_{L^\infty(Q_4)}$ and
$(H^*)^s\phi\in L^1$. To verify the latter, write
$S_r^*\phi(x,t)=e^{r\Delta_x}\phi(x,t+r)$, with all $L^1$ norms
below taken over $\R^{n+1}$. The heat kernel has mass one and time
translation preserves the $L^1$ norm, so
$\|S_r^*g\|_1\le\|g\|_1$. Consequently,
\[
 \|\phi-S_r^*\phi\|_1\le2\|\phi\|_1.
\]
On the other hand, $\frac{d}{dr}S_r^*\phi=-S_r^*(H^*\phi)$ gives
\[
 \phi-S_r^*\phi=\int_0^rS_q^*(H^*\phi)\,dq,
 \qquad
 \|\phi-S_r^*\phi\|_1\le r\|H^*\phi\|_1.
\]
Using this bound for $0<r<1$ and the preceding bound for $r\ge1$,
the adjoint semigroup formula yields
\[
 \begin{aligned}
 \|(H^*)^s\phi\|_1
 &\le\frac1{|\Gamma(-s)|}
   \int_0^\infty\|\phi-S_r^*\phi\|_1\frac{dr}{r^{1+s}}\\
 &\le\frac1{|\Gamma(-s)|}
   \left(\|H^*\phi\|_1\int_0^1r^{-s}\,dr
         +2\|\phi\|_1\int_1^\infty r^{-1-s}\,dr\right)\\
 &=\frac1{|\Gamma(-s)|}
   \left(\frac{\|H^*\phi\|_1}{1-s}
         +\frac{2\|\phi\|_1}{s}\right)<\infty,
 \end{aligned}
\]
since $0<s<1$. As $\phi$ is fixed, this proves the bound used
in \eqref{eq:tail-average}.

We claim that, for every multi-index $\beta$ and integer $k\ge0$,
\begin{equation}\label{eq:kernel-domination}
 |\partial_x^\beta\partial_t^kK(x-z,t-\tau)|
 \le C_{\beta,k}K(x'-z,t'-\tau)
\end{equation}
whenever $(x,t)\in Q_{5/2}$, $(x',t')\in A$, $\tau<t$, and
$h(z,\tau)\ne0$. Put
$\xi=z-x$, $b=x'-x$, $r=t-\tau$, and $\delta=t'-t$.
Then $3/4<\delta<57/4$, $|b|\le7/2$, and completing the square gives
\begin{equation}\label{eq:kernel-ratio}
 \frac{K(\xi,r)}{K(\xi-b,r+\delta)}
 =\left(1+\frac\delta r\right)^\nu
 \exp\!\left(-\frac{\delta|\xi+(r/\delta)b|^2}{4r(r+\delta)}
                  +\frac{|b|^2}{4\delta}\right).
\end{equation}
Moreover, direct differentiation yields, with $\ell=|\beta|+2k$,
\begin{equation}\label{eq:kernel-derivatives}
 |\partial_\xi^\beta\partial_r^kK(\xi,r)|
 \le C_{\beta,k}r^{-\ell/2}
       (1+|\xi|/\sqrt r)^\ell K(\xi,r).
\end{equation}
For $r\ge r_0:=1/100$, set
$Y=(\xi+(r/\delta)b)/\sqrt{r(r+\delta)}$. Then
\[
 r^{-\ell/2}(1+|\xi|/\sqrt r)^\ell
 =(r^{-1/2}+|\xi|/r)^\ell\le C(1+|Y|)^\ell,
\]
because $|\xi|/r\le\sqrt{(r+\delta)/r}\,|Y|+|b|/\delta$.
Thus \eqref{eq:kernel-ratio}--\eqref{eq:kernel-derivatives} bound the
differentiated quotient by $C(1+|Y|)^\ell e^{-c|Y|^2}\le C$.
For $0<r<r_0$, we have $\tau\in(-9,9)$, so $h(z,\tau)\ne0$
implies $|z|\ge3$ and $|\xi|\ge1/2$. Since
$|(r/\delta)b|\le1/4$, we have $|a+(r/\delta)b|\ge|a|/2$.
The same quotient is therefore bounded by
\[
 C r^{-\nu-\ell/2}(1+|a|/\sqrt r)^\ell
 e^{-c|a|^2/r}\le C r^{-N}e^{-c'/r}\le C.
\]
This proves \eqref{eq:kernel-domination}. 

Extend $K$ by zero for $r\le0$. Support separation removes the
singularity at $r=0$, and \eqref{eq:kernel-domination} justifies
differentiation under the tail integral. Since $h\ge0$ and $t<t'$, it gives
\[
 |\partial_x^\beta\partial_t^kT(x,t)|
 \le C_{\beta,k}\int_{-\infty}^t\!\int_{\R^n}
 K(x'-z,t'-\tau)h(z,\tau)\,dz\,d\tau
 \le C_{\beta,k}T(x',t').
\]
Averaging over $(x',t')$ against $\phi$ and using
\eqref{eq:tail-average}, we obtain
\begin{equation}\label{eq:smooth-tail}
 \|\partial_x^\beta\partial_t^kT\|_{L^\infty(Q_{5/2})}
 \le C_{\beta,k}L.
\end{equation}
In particular, $T$ is smooth there and
$\|T\|_{C^{\eta,\eta/2}(Q_{5/2})}\le C_\eta L$ for $0<\eta<1$.

We now apply \cite[Theorem~1.20]{ST2017}, whose hypotheses require global
boundedness and membership in the $L^2$ operator domain of $H^s$.
To justify its use for our distributional solution, take a nonnegative, compactly supported
smooth space-time approximate identity and set
$w_\epsilon=\rho_\epsilon*w$. Then
$w_\epsilon\in C_c^\infty\subset\operatorname{Dom}(H^s)$,
$\|w_\epsilon\|_\infty\le\|w\|_\infty$, and translation invariance gives
\[
 H^s w_\epsilon=(f+T)*\rho_\epsilon
 \quad\text{in }Q_{12/5}
\]
for sufficiently small $\epsilon$. The source norms on this cylinder
are uniformly controlled by those of $f+T$ on $Q_{5/2}$.
The interior estimates, followed by compactness and passage to the
limit, therefore give
\[
 \|u\|_{C^{\theta,\theta/2}(Q_2)}
 =\|w\|_{C^{\theta,\theta/2}(Q_2)}
 \le C\bigl(\|w\|_\infty+\|f+T\|_{L^\infty(Q_{5/2})}\bigr)
 \le CL
\]
for $0<\theta<\min\{2s,1\}$. At $2s=1$, the H\"older--Zygmund
estimate implies every such bound with $\theta<1$.
If $f\in C^{\eta,\eta/2}(Q_4)$ and $\lambda=2s+\eta\in(0,2)\setminus\{1\}$,
the Schauder estimate similarly yields
\[
 \|u\|_{C^{\lambda,\lambda/2}(Q_2)}
 \le C\bigl(\|w\|_\infty+
       \|f+T\|_{C^{\eta,\eta/2}(Q_{5/2})}\bigr)
 \le C\bigl(\|u\|_{L^\infty(Q_4)}+
       \|f\|_{C^{\eta,\eta/2}(Q_4)}\bigr).
\]
For $1<\lambda<2$, apply the separate spatial and temporal estimates
to $w_\epsilon$ on $Q_{9/4}$ and write
$\mathcal N=\|u\|_{L^\infty(Q_4)}+\|f\|_{C^{\eta,\eta/2}(Q_4)}$.
Comparing derivatives with spatial difference quotients at scale
$r=|t-\tau|^{1/2}$ gives, for nearby times and $x\in B_2$,
\[
 |\nabla_xw_\epsilon(x,t)-\nabla_xw_\epsilon(x,\tau)|
 \le C\mathcal N\bigl(r^{\lambda-1}
              +r^{-1}|t-\tau|^{\lambda/2}\bigr)
 =C\mathcal N|t-\tau|^{(\lambda-1)/2}.
\]
Boundedness of the gradients handles larger time differences, completing
also the temporal H\"older estimate for $\nabla_xw$ after passage to the limit.
\end{proof}

We now combine this local regularity estimate with doubling and the
bounded-solution Liouville theorem to prove Theorem~\ref{main-thm-2}. 

\begin{proof}[Proof of Theorem~\ref{main-thm-2}]
Set $\alpha=2s/(p-1)$. If \eqref{estimate-universe} fails,
there are domains $\Omega_j$, solutions $u_j$, and points
$P_j\in\Omega_j$ such that, writing
$m_j=u_j^{1/\alpha}$ and $d_j(P)=\operatorname{dist}_p(P,\Omega_j^c)$,
\[
    m_j(P_j)d_j(P_j)>4j.
\]
Here is the doubling construction underlying \cite{PQS12007,PQS2007}.
Start with $P_{j,0}=P_j$. At step $k$, stop if
\[
    \sup_{Q_{j/m_j(P_{j,k})}(P_{j,k})}m_j\le2m_j(P_{j,k});
\]
otherwise choose $P_{j,k+1}$ in that cylinder with
$m_j(P_{j,k+1})>2m_j(P_{j,k})$.
Inductively $m_j(P_{j,k})\ge2^k m_j(P_j)$, and every point $Q$
in the cylinder tested at step $k$ satisfies
\[
 d_p(Q,P_j)<\sum_{\ell=0}^k\frac{j}{m_j(P_{j,\ell})}
 \le\frac{j}{m_j(P_j)}\sum_{\ell=0}^k2^{-\ell}
 <\frac{2j}{m_j(P_j)}<\frac{d_j(P_j)}2.
\]
Thus all choices and test cylinders lie in the fixed compact set
$\overline Q_{2j/m_j(P_j)}(P_j)\Subset\Omega_j$.
An infinite construction would contradict local boundedness of $m_j$.
At the stopping point $Z_j=(z_j,\tau_j)$, therefore,
\begin{equation}\label{eq:doubling}
 M_j:=m_j(Z_j)\ge m_j(P_j),\qquad
 Q_{j/M_j}(Z_j)\subset\Omega_j,\qquad
 \sup_{Q_{j/M_j}(Z_j)}m_j\le2M_j.
\end{equation}

Put $r_j=M_j^{-1}$ and define, on all of $\R^n\times\R$,
\[
    v_j(x,t)=r_j^\alpha u_j(z_j+r_jx,\tau_j+r_j^2t).
\]
Since $\alpha+2s=\alpha p$, \eqref{eq:doubling} and scaling give
\[
    H^s v_j=v_j^p\quad\text{in }Q_j,
    \qquad v_j(0,0)=1,
    \qquad 0\le v_j\le B:=2^\alpha\quad\text{in }Q_j.
\]
Translations and parabolic dilations preserve $\mathcal L_s$, so $v_j\in\mathcal L_s$;
also $v_j\ge0$ globally. 

Fix $0<\theta<\min\{2s,1\}$. On every $Q_4(P)\subset Q_j$,
Proposition~\ref{prop:local-regularity} gives
\[
 \|v_j\|_{C^{\theta,\theta/2}(Q_2(P))}\le C(B+B^p).
\]
Since $a\mapsto a^p$ is Lipschitz on $[0,B]$, the same local
H\"older bound holds for $v_j^p$.
Choose $0<\eta<\min\{\theta,2-2s\}$ with $2s+\eta\ne1$,
and put $\gamma=2s+\eta$.
Applying the Schauder estimate on $Q_2(P)\supset Q_1(P)$ gives
\[
 \|v_j\|_{C^{\gamma,\gamma/2}(Q_1(P))}\le A,
\]
where $A$ depends only on $n,s,p$, and is independent of $j$ and $P$.
Arzel\`a-Ascoli theorem and a diagonal extraction give, for some
$0<\varepsilon<\eta$ with $2s+\varepsilon\ne1$,
\[
 v_j\longrightarrow v\quad\text{in }
 C_{\mathrm{loc}}^{2s+\varepsilon,s+\varepsilon/2},\qquad
 0\le v\le B,\qquad v(0,0)=1.
\]
For each fixed $P$, $Q_4(P)\subset Q_j$ eventually. Passing the
same bound to the limit therefore gives
\begin{equation}\label{eq:uniform-limit-holder}
 \sup_{P\in\R^n\times\R}
 \|v\|_{C^{\gamma,\gamma/2}(Q_1(P))}\le A.
\end{equation}
In particular, $v\in\mathcal L_s$: for $\nu=n/2+1+s$, boundedness gives
\begin{equation}\label{eq:limit-tail}
 \int_{-\infty}^t\!\int_{\R^n}
 \frac{|v(z,\tau)|e^{-|x-z|^2/[4(t-\tau)]}}
 {1+(t-\tau)^\nu}\,dz\,d\tau
 \le (4\pi)^{n/2}B\int_0^\infty\frac{r^{n/2}}{1+r^\nu}\,dr<\infty.
\end{equation}

To bound $H^sv$, write
$\delta_z^2v(x,t)=v(x+z,t)+v(x-z,t)-2v(x,t)$.
The evenness of $G_r$ gives
\[
 v(x,t)-S_rv(x,t)
 =v(x,t)-v(x,t-r)
   -\frac12\int_{\R^n}G_r(z)\delta_z^2v(x,t-r)\,dz.
\]
For $0<r<1$, \eqref{eq:uniform-limit-holder} bounds the time
difference by $Ar^{\gamma/2}$.
For $|z|<1$, it gives $|\delta_z^2v|\le CA|z|^\gamma$:
if $\gamma<1$, use the spatial H\"older estimate; if $\gamma>1$, use
\[
 \delta_z^2v(x,t)
 =\int_0^1[\nabla_xv(x+qz,t)-\nabla_xv(x-qz,t)]\cdot z\,dq
\]
and the $(\gamma-1)$-H\"older continuity of $\nabla_xv$.
For $|z|\ge1$, boundedness gives
$|\delta_z^2v|\le4B\le4B|z|^\gamma$.
Since $\int G_r(z)|z|^\gamma\,dz=C_{n,\gamma}r^{\gamma/2}$,
\[
 \|v-S_rv\|_\infty\le C(A+B)r^{\gamma/2}\quad(0<r<1),
 \qquad \|v-S_rv\|_\infty\le2B\quad(r\ge1).
\]
The semigroup formula now yields
\begin{equation}\label{eq:limit-operator-bound}
 \begin{aligned}
 \|H^sv\|_\infty
 &\le\frac1{|\Gamma(-s)|}
 \left(C(A+B)\int_0^1r^{\gamma/2-1-s}\,dr
       +2B\int_1^\infty r^{-1-s}\,dr\right)\\
 &\le\frac1{|\Gamma(-s)|}
 \left(\frac{C(A+B)}{\gamma/2-s}+\frac{2B}{s}\right)<\infty.
 \end{aligned}
\end{equation}
Here $\gamma>2s$ makes the first integral finite.

For each fixed $P$, the equation for $v_j$ holds at $P$ for all
sufficiently large $j$. Consequently,
\[
 \lim_{j\to\infty}H^s v_j(P)=v(P)^p,\qquad
 \liminf_{j\to\infty}H^s(v_j-v)(P)
 =v(P)^p-H^s v(P)\ge-\|H^s v\|_\infty.
\]
Recall $v$ is bounded, we may therefore apply \cite[Theorem~1.1 and Remark~1.3]{CGLO2026}, obtaining
\[
    H^s v=v^p+b\quad\text{in }\R^n\times\R
\]
for some constant $b\ge0$.

Applying Lemma~\ref{lem:representation} to $v$ with
$f=v^p+b$ gives, $b=0$.
Hence $v$ is a bounded
nonnegative entire solution of $H^s v=v^p$.
Proposition~\ref{prop:bounded-zero} implies $v\equiv0$, contradicting $v(0,0)=1$.
\end{proof}

\begin{proof}[Proof of Theorem~\ref{main-thm-1}]
Fix $P\in\R^n\times\R$. Applying Theorem~\ref{main-thm-2}
in $Q_R(P)$ gives $u(P)\le CR^{-2s/(p-1)}$.
Letting $R\to\infty$ proves $u(P)=0$.
\end{proof}

\medskip
{\textbf{AI disclosure statement.}}
The authors initially established the main results of the present paper by deriving an energy monotonicity formula for the extension of solutions and by employing the method in Quittner's article \cite{Quittner2021}. However, the argument was excessively intricate. Through iterative interactions with ChatGPT (OpenAI), the concise and simplified proof presented here was developed. All mathematical results and
proofs were reviewed and verified by the authors, who take full responsibility for the content of the manuscript.

\medskip
{\textbf{Acknowledgments.}} 
Guo is partially supported by the National Natural Science Foundation of China (Grant No. 12501145), the Natural Science Foundation of Shanghai (No. 25ZR1402207),   the China Postdoctoral Science Foundation (No. 2025T180838 and No. 2025M773061), the Postdoctoral Fellowship Program of CPSF (No. GZC20252004), and the Institute of Modern Analysis-A Frontier Research Center of Shanghai.
Ouyang and Xie are partially supported by the National Natural Science Foundation of China (Grant No. W2531006, 12250710674 and 12031012) and the Institute of Modern Analysis-A Frontier Research Center of Shanghai. 
\medskip

{\textbf{Data availability statement.}} Data will be made available on reasonable request.
\medskip

{\textbf{Conflict of interest statement.}} There is no conflict of interest.

\end{document}